\documentclass[12pt]{article}
\usepackage{amsfonts, amssymb,amsmath, amsthm, dsfont, url,epsfig,fancyhdr,graphicx,mathrsfs,eucal}
\usepackage{array}
\usepackage{eucal}
\newcommand\comment[1]{}                %  Silent version.
\renewcommand\comment[1]{\emph{[#1]}}           %  Comment revealed.
\def\int{\mathop{\rm int}}

\newtheorem{theorem}{Theorem}

\newtheorem{lemma}{Lemma}[section]

\newtheorem{proposition}[theorem]{Proposition}

\theoremstyle{remark}

\theoremstyle{definition}
\newtheorem{definition}[theorem]{Definition}

\renewcommand{\phi}{\varphi}          % Personal preferences.
\renewcommand{\epsilon}{\varepsilon}

\newcommand\sgn{\operatorname{sgn}}

\renewcommand\ker{\operatorname{\mathfrak{Ker}}}
\newcommand\conv{\operatorname{conv}}
\newcommand\aff{\operatorname{aff}}

\renewcommand\vert{\operatorname{\mathsf{vert}}}
\newcommand\vertex{\operatorname{\mathsf{vert}}}

\newcommand\GCD{\operatorname{GCD}}

\newcommand\Sym{\operatorname{\mathtt{Sym}}}

\renewcommand\a{\mathbf{a}}
\renewcommand\b{\mathbf{b}}
\newcommand\e{\mathbf{e}}
\renewcommand\j{\mathbf{j}}

\newcommand\uu{\mathbf{u}}  % "\u" is preassigned another meaning on Latex
\newcommand\vv{\mathbf{v}} % "\v" is preassigned another meaning on Latex
\newcommand\cc{\mathbf{c}}

\newcommand\x{\mathbf{x}}
\newcommand\y{\mathbf{y}}
\newcommand\z{\mathbf{z}}

\newcommand\sV{\mathsf{V}}

\newcommand\bV{\mathbf{V}}

\newcommand\cP{\mathcal{P}}

\newcommand\R{\mathbb{R}}    %reals
\newcommand\N{\mathbb{N}}    %naturals
\newcommand\Q{\mathbb{Q}}   %rationals
\newcommand\Z{\mathbb{Z}}    %integers
\newcommand\E{\mathbb{E}}     %Euclidean space, i.e. real affine space with fixed positive quadratic metric
\newcommand\LL{\mathbb{L}}
\newcommand\Mu{\mathrm{M}}

\begin{document}

\author{\large Robert Erdahl, Andrei Ordine, and Konstantin Rybnikov\footnote{corresponding author}}
\title{\large Perfect Delaunay Polytopes and Perfect Quadratic Functions on Lattices}
\date{\today}
\maketitle
\begin{abstract}

A polytope $D$ whose vertices belong to a lattice of rank $d$ is \textit{Delaunay} if there is 
a circumscribing $d$-dimensional ellipsoid, $\mathcal{E}$, with interior free of lattice points 
so that the vertices of $D$ lie on $\mathcal{E}$.  If in addition, the ellipsoid $\mathcal{E}$ 
is uniquely determined by $D$, we call $D$ \emph{perfect}.  That is, a perfect Delaunay 
polytope is a lattice polytope with a circumscribing empty ellipsoid $\mathcal{E}$, where the 
quadratic surface $\partial \mathcal{E}$ both contains the vertices of $D$ and is determined by 
them.  We have been able to construct infinite sequences of perfect Delaunay polytopes, one 
perfect polytope in each successive dimension starting at some initial dimension; we have been 
able to construct an infinite number of such infinite sequences. Perfect Delaunay polytopes 
play an important role in the theory of Delaunay polytopes, and in Voronoi's theory of lattice 
types.
\end{abstract}

\section{ \large Introduction\label{introduction}}
Consider the lattice $\mathbb{Z}^{d}$, and a lattice polytope $D$. Then, $D$ is a 
\emph{Delaunay polytope} if it can be circumscribed by an \emph{ellipsoid} $\mathcal{E}$ so 
that: (1) $\mathcal{E}$ has no interior $\mathbb{Z}^{d}$-elements, and (2) the only 
$\mathbb{Z}^{d}$-elements on $\mathcal{E}$ are the vertices of $D$. We will say that the 
circumscribing ellipsoid $\mathcal{E}$, with these two properties, is an \emph{empty 
ellipsoid}. If $D$ is a Delaunay polytope with circumscribing empty ellipsoid $\mathcal{E}$, 
then $D$ is \emph{perfect} if the quadratic surface $ \partial \mathcal{E}$ is uniquely determined by the vertex 
set of $D$.

If $D$ is Delaunay, with empty ellipsoid $\mathcal{E}$, the vertices of $D$ are given by 
$\partial \mathcal{E} \bigcap D = \partial \mathcal{E} \bigcap \mathbb{Z}^{d}$.  Typically, the 
quadratic surface $\partial \mathcal{E}$, through this vertex set, can be deformed 
continuously.  It is only in the case where $D$ is a perfect Delaunay polytope that this 
surface is uniquely determined by $\partial \mathcal{E} \bigcap D$.  We refer to $\partial 
\mathcal{E}$ as a \emph{perfect ellipsoid}, and by extension, to refer to $D$ as a\emph{ perfect 
Delaunay polytope. } 

We have studied perfect Delaunay polytopes by constructing infinite
sequences of them, one perfect Delaunay polytope in each successive
dimension, starting at some initial dimension.  We have been able to
construct an infinite number of infinite sequences of perfect
Delaunay polytopes.

The perfect Delaunay polytopes in each sequence have similar combinatorial and geometric 
properties.  One of our constructions is the sequence of  $G_6$-\emph{topes,} $G_{6}^d,\;d=6,7,...$ , 
with the initial term the semiregular 6-dimensional Gossett polytope $G_6$ with 27 vertices 
($2_{21}$ in Coxeter's notation); each $G_6$-tope is asymmetric with respect to central 
inversion, and $G_{6}^d$ has $\frac{d(d+1)}{2}+d$ vertices.  Each term in the sequence is a 
combinatorial analogue of $G_{6}=2_{21}$. Another of our constructions is the sequence of 
\emph{$G_7$-topes}, $G_{7}^d,\;d=7,8,...$, with the initial term the 7-dimensional Gossett polytope 
$G_7$ with 56 vertices ($3_{31}$ in Coxeter's notation); each $G_7$-tope is centrally-symmetric, 
and $G_{7}^d$ has $2 \binom{d+1}{2}$ vertices. Each term is again a combinatorial analogue of the 
initial Gosset polytope $G_7=3_{31}$.  Just as the 6-dimensional Gossett polytope can be 
represented as a section of the 7-dimensional one, each $G^{d}_6$ of the asymmetric sequence can 
be represented as a section of the $G^{d+1}_7$ of the symmetric sequence. 

\medskip

\noindent \textbf{Properties and Connections:} Perfect Delaunay polytopes are fascinating 
geometric objects. There is a number of other reasons why we  studied these polytopes.

\medskip

\begin{enumerate}

\item Perfect Delaunay polytopes are maximal Delaunay polytopes  -- it follows immediately from the
definition that a perfect Delaunay polytope  cannot fit properly inside another Delaunay 
polytope of the same dimension.

\item Each perfect Delaunay polytope $D$ uniquely determines a positive definite quadratic form
on $\mathbb{Z}^{d}$, which in turn, determines a Delaunay tiling of $\mathbb{Z}^{d}\otimes \R$ 
with $D$ as one of its  Delaunay tiles (details on how these tilings are constructed are given 
in the following section). Therefore, each of our sequences of perfect Delaunay polytopes 
determines a corresponding sequence of Delaunay tilings. Just as the perfect Delaunay polytopes 
in the sequence are combinatorial analogues, the corresponding Delaunay tilings have somewhat 
similar combinatorial properties.

Each perfect Delaunay polytope in the sequence of $G_6$-topes uniquely determines a Delaunay 
tiling that is an analogue of that for the root lattice $E_{6}$ in $\E^6$.  Similarly, each 
term in the sequence of $G_7$-topes uniquely determines a Delaunay tiling that is an analogue of 
that for the root lattice $E_{7}$  in $\E^7$.

\item A  quadratic form $\varphi$ in $d$ indeterminates is called perfect (after Voronoi) if it is positive definite and is uniquely
determined by its minimal vectors through the linear system of equations $\{\chi[\mathbf{m}]=\min_{\neq0} 
\phi\;\;\vline\; \;\mathbf{m}\: \text{is a minimal vector for} \:\varphi\}$, where the coefficients of $\chi$ are the variables, $\min_{\neq 0} \phi$  is the minimal value 
achieved by $\varphi$ on the non-zero elements of $\mathbb{Z}^{d}$, and each equation corresponds a minimal vector $\mathbf{m} \in \Z^d$ for 
$\varphi$. If $\varphi$ is perfect, then the only solution is 
$\chi=\varphi$.

There is an inhomogeneous analogue to perfect form that plays an important role in the study of 
Delaunay polytopes.  Let $\mathcal{E}$ be an empty ellipsoid around a Delaunay polytope $D$, 
defined by inequality $f(\x)\le 0$.  Since $\mathcal{E}$ is empty, $f$ is non-negative on 
$\mathbb{Z}^{d}$, and assumes its minimum value zero on $\partial \mathcal{E} \bigcap 
\mathbb{Z}^{d}$.  If $D$ is perfect then the coefficients of the inhomogeneous quadratic 
function $f$ are determined, up to a scale factor, by the linear system of equations 
$\{f(\mathbf{z})=0\;\;\vline\;\; \mathbf{z}\in \mathcal{E} \bigcap \mathbb{Z}^{d}\}$.  The rank of this linear 
system must then be $\frac{d(d+1)}{2}+d$, which requires a perfect Delaunay polytope to have at 
least this number of vertices.  The vertices of a perfect Delaunay polytope are used to 
determine the coefficients of $f$ in much the same way that the minimal vectors are used to 
determine the coefficients of a perfect form. Therefore, when $D$ is perfect it is natural to 
refer to $f$ as a \emph{perfect quadratic function}, and to refer to $\mathcal{E}$ as a 
\emph{perfect ellipsoid}. Clearly $\mathbf{0}$ is perfect Delaunay in $\Z^0$ and 
$\mathbf{[0,1]}$ is perfect Delaunay in $\Z^1$. Next perfect Delaunay polytope lives in 
dimension 6.

The conditions that a perfect quadratic function must satisfy are more severe than those for a 
perfect form, which requires at least $\frac{n(n+1)}{2}$ pairs of minimal vectors.  This leads 
to the speculation that the growth of numbers of arithmetic types of perfect quadratic 
functions with dimension will be much slower than for numbers of types of perfect forms.  This 
has been borne out by the data now available on perfect Delaunay polytopes (see 
\url{http://www.liga.ens.fr/~dutour/}).  In fact, the growth of numbers of arithmetic types of 
perfect function seems much slower that that for most classes of geometric objects studied in 
geometry of numbers. 

Perfect forms have been introduced by Korkine and Zolotareff (1873), who classified all the 
perfect forms through dimension five (the term \verb"perfect" was coined by Voronoi). In 
addition, they found several perfect forms in higher dimensions, including the forms 
corresponding to the root lattices $A_n$, $D_n$ (both for $n \in \N$), $E_{6},E_{7},$ and $E_{8}$, although they used completely
different roman letters to denote these lattices.  Conway and Sloane (1988) give a 
survey of perfect forms through dimension seven. More recently, the study of perfect forms has 
been taken up by Martinet who published a comprehensive monograph on the topic in 2003. As of 
the end of 2006 perfect forms have been classified in dimensions through eight (see 
\url{http://fma2.math.uni-magdeburg.de/~latgeo/}).

\item A significant role for perfect Delaunay polytopes is
to provide geometric labels for a class of \emph{edge forms} for the \emph{homogeneous domains} 
of Delaunay tilings of $\mathbb{Z}^{d}\otimes \R$.  These domains, also referred to as L-type 
domains, are polyhedral cones of positive forms, and the edge forms are those that lie on the 
extreme rays. L-type domains have been first introduced and studied  by  Voronoi, and they play 
a central role in his classification theory for lattices, his theory of lattice types. Delaunay  
associated these polyhedral cones with Delaunay tilings. Prior to the discovery of the sequence 
$\{G_6^{d}\}$ of $G_6$-topes by Erdahl and Rybnikov in 2001 only finitely many arithmetic types of 
edge forms of \emph{full rank} were known. All of the infinite sequences of perfect Delaunay polytopes 
that we report on correspond to infinite sequences of edge forms of full rank for L-type 
domains. Later infinite series of perfect Delaunay polytopes were constructed by Dutour (2005) and  Grishukhin (2006). Our multiparametric series overlap with the 1-parametric series of Dutour only in small dimensions (6,7,8). Our series contain all of the series of Grishukhin, although his approach is different. The remarkable property of Dutour's series is that the size of the vertex set of Dutour's polytope grows exponentially with the dimension, for his $d$-polytope $ED_d$, for even $d$'s, contains a section isometric to the $(d-1)$-halfcube (for odd $d$'s it has 2 sections isometric to the $(d-2)$-halfcube). Dutour conjectured that this polytope $ED_d$ has the largest number of vertices among all perfect Delaunay polytopes of dimension  $d$.

Perfect Delaunay polytopes, and the infinite cylinders with perfect
Delaunay polytopes as bases, also provide labels for \emph{edge
functions} for the \emph{inhomogeneous domains} of Delaunay
polytopes.  Associated with each lattice Delaunay polytope is a polyhedral
cone of inhomogeneous quadratic functions, and we refer to the
functions lying on extreme rays as edge functions.  These domains
play an important role in the structure theory of Delaunay polytopes.

\item Delaunay tilings of lattices and Voronoi's theory of L-types  are intrinsically connected to the theory of moduli spaces of abelian varieties. On a philosophical level this connection stems from the idea of parametric approach to the study of geometric objects. In geometry of numbers and number theory this approach was pioneered by Hermite (1850). Korkine, Zolotareff, and Voronoi  referred to this approach as  Hermite's method of continuous variation of parameters. It is worth noting  that Voronoi's interest in geometry of numbers developed as a result of his work on irrationalities of the third degrees and, in particular, elliptic curves. The curous reader may consult Delone [Delaunay] and Faddeev (1964). For recently discovered connections between moduli spaces of abelian varieties and Voronoi's theories of perfect and  of L-types see Alexeev (2002) and Shepard-Barron (2005). We are not aware for any results interpreting perfect perfect Delaunay polytopes and L-types defined by them, but we certainly anticipate them.

The following section is devoted to the details of how perfect
Delaunay polytopes relate to these fundamental polyhedral cones in
geometry of numbers, the homogeneous domain of a Delaunay tiling,
and the inhomogeneous domain of a Delaunay polytope.

\end{enumerate}

\noindent \textbf{Historical notes:} Perfect Delaunay polytopes were first considered by Erdahl 
(1975) in connection with lattice polytopes arising from the quantum mechanics of many 
electrons.  He observed (1975) that Delaunay tilings of 0- and 1-dimensional lattices consist 
entirely of perfect Delaunay polytopes, and showed that the Gosset polytope $G_6=2_{21}$ with 
27 vertices was a perfect Delaunay polytope in the root lattice $E_{6}$. \ He also showed that 
there were no perfect Delaunay polytopes in dimensions 2, 3, and 4. \ These results were 
further extended by Erdahl (1992) by showing that the 7-dimensional Gosset polytope $G_7=3_{31}$ with 
56 vertices is perfect, and that there are no perfect Delaunay polytopes of 
dimension $d$ for $1<d<6$. Erdahl also proved that $G_6$ and $G_7$ are the only perfect
Delaunay polytopes in the Delaunay tilings for the root lattices. Deza, Grishukhin, and Laurent 
(1992-1997, in various combinations) found more examples of perfect Delaunay polytopes in dimensions 15, 16, 22, and 23, but all of those seemed to be sporadic.

The first construction of infinite sequences of perfect Delaunay polytopes was described at the 
Conference dedicated to the Seventieth Birthday of Sergei Ryshkov (Erdahl, 2001), and later 
reported by Rybnikov (2001) and Erdahl and Rybnikov (2002). Perfect Delaunay polytopes have 
been classified up to dimension 7 -- Dutour (2004) proved that $G_6=2_{21}$ is the only perfect 
polytope for $d=6$. It is strongly suspected that the existing lists of seven and eight 
dimensional perfect Delaunay polytopes are complete (see \url{http://www.liga.ens.fr/~dutour}).

\section{\large Quadratic Polynomials on Lattices\label{forms}}

 Perfect Delaunay polytopes play an important role in the
\emph{ L-type reduction theory of Voronoi and Delaunay} for positive (homogeneous) quadratic 
forms. An L-type domain is the collection of all possible positive quadratic  forms that give the 
same Delaunay tiling $\mathcal{D}$ for $\mathbb{Z}^{d}$. L-type domains are relatively open 
polyhedral cones, with boundary cells that are also L-type domains -- these conical cells fit 
together to tile the cone of positive quadratic forms, which is described in the next section 
in  detail. Simplicial Delaunay tilings label the full dimensional conic "tiles'', and all 
other possible Delaunay tilings label the lower dimensional cones. \

The significance of extreme L-types is much due to their relation to the structure of Delaunay 
and Voronoi tilings for lattices. The Delaunay tilings that correspond to edge forms play an 
important role: \ \emph{The Delaunay tiling for an L-type domain is the intersection of the 
Delaunay tilings for edge forms for the L-type} (Erdahl, 2000).  There is a corresponding dual 
statement on the
structure of Voronoi polytopes: \emph{The Voronoi polytope }$V_{\varphi }$%
\emph{\ for a form }$\varphi $\emph{\ contained in an L-type domain $\mathcal{L}$ is a weighted Minkowski sum 
of linear transforms of Voronoi polytopes for each of the edge forms of $\mathcal{L}$. }\  The latter dual 
result was first established by H.-F. Loesch  in  his 1990 doctoral dissertation, although it 
was first published by Ryshkov (1998, 1999), who independently rediscovered Loesch's theorem; 
this dual result was given a shorter and simpler proof by Erdahl (2000). \

Edge forms that are interior to the cone of positive forms are rare in low dimensions. They first 
occur in dimension 4: $D_{4}$ has an extreme $L $-type. A good proportion, but not all, of the 
edge forms appearing in lower dimensions relate either directly or indirectly to perfect 
Delaunay polytopes. As shown by Dutour and Vallentin (2005) this situation does not persist in 
higher dimensions: there is an ``explosion" of the number of inequivalent interior edge forms in six dimensions, and only a tiny 
fraction of these are inherited from perfect Delaunay polytopes.

The Voronoi and Delaunay tilings for point lattices are constructed using the Euclidean metric, 
but are most effectively studied by injectively mapping the lattice into $\mathbb{Z}^{d}$, and 
replacing the Euclidean metric by an
equivalent metrical form. \ \ For a $d$-dimensional point lattice $
\Lambda $ with a basis $\mathbf{b}_{1},\mathbf{b}_{2},...,\mathbf{b}_{d}$ this is done as 
follows. \ \ A lattice vector with coordinates $z_{1},z_{2},...,z_{d}$ relative to this basis 
can be written as $\mathbf{z}=\mathbf{Bz}$, where 
$\mathbf{B=[\mathbf{b}_{1},\mathbf{b}_{2},...,\mathbf{b}_{d}]}$ is the basis matrix
and $\mathbf{z}$ is the column vector given by $\mathbf{[}%
z_{1},z_{2},...,z_{d}]^{t}$. The squared Euclidean length is given by $|\mathbf{v }|^{2}=%
\mathbf{z}^{t}\mathbf{B}^{t}\mathbf{Bz=}\varphi _{\mathbf{B}}(\mathbf{z})$. \ This squared 
length can equally well be interpreted as the squared length of the integer vector 
$\mathbf{z}\in \mathbb{Z}^{d}$ relative to the form $\varphi _{\mathbf{B}}$. \ Therefore, the 
Dirichlet-Voronoi and Delaunay tilings for $\Lambda $, constructed using the Euclidean 
metric, can be studied using the corresponding Dirichlet-Voronoi and Delaunay tilings for 
$\mathbb{Z}^{d}$ constructed using the form $\varphi _{\mathbf{B}}$ as metric. \ \ Moreover,
variation of the Dirichlet-Voronoi and Delaunay tilings for $\Lambda$ in response to 
variation of the lattice basis  can be studied by varying the metric  $\varphi $
for the fixed lattice $\mathbb{Z}^{d}$. In this section we will keep the lattice fixed at 
$\mathbb{Z}^{d}$, and vary the quadratic form $\varphi $. We call a form $\varphi$ positive if  
$\varphi[\z] \ge 0$ for any $\z$, and we call it positive  definite if $\varphi[\z]>0$ for any 
$\z\neq \mathbf{0}$. Same terminology is applied to all numbers and functions. Positive, but 
not positive definite forms are referred to as positive semidefinite, or just semidefinite.

\noindent \textbf{The inhomogeneous domain of a Delaunay polytope: }   Let $\cP^d$ be the cone 
of real quadratic polynomials defined by:
\begin{equation*}
\cP^d = \{\, f \in \R [x_1,\dots,x_d] \;\; \vline \;\; \deg f = 2, \;\; \forall \mathbf{z} \in 
\mathbb{Z}^{d}: f(\mathbf{z}) \geq 0 \, \}.
\end{equation*}
The condition $\forall \mathbf{z}\in \mathbb{Z}^{d}: f(\mathbf{z})\geq 0$ requires  the 
quadratic part of $f$ to be positive, and requires  any subset of $\mathbb{R}^{d}$ where $f$ 
assumes negative values to be free of $\mathbb{Z}^{d}$-elements. The real quadric determined by 
the equation $f(\mathbf{x})=0$, where $f \in \cP^d$, might be empty set  or it might have the form
\begin{equation*}
\partial \mathcal{E}_{f}=\partial \mathcal{E}_{0}\times K,
\end{equation*}
where $\mathcal{E}_{0}$ is an ellipsoid of dimension $n\ge 0$ and $K$ a complementary subspace, of dimension $d-n$. \ The latter case is 
the interesting one -- depending on the dimension of $K$, $ \mathcal{E}_{f}$ is either an empty 
(of lattice points) ellipsoid or an empty cylinder with ellipsoidal base. 

For any function $f \in \cP^d$ we denote by $\sV(f)$ the set $\{\ \z \in \Z^d \; \vline\; 
f(\z)=0  \}$. In the case where the surface $\partial \mathcal{E}_{f}=\{\ \x \in \Z^d \otimes 
\R \; \vline\; f(\x)=0\}$ includes integer points and is bounded, $\sV(f)=\mathcal{E}_{f}\, 
\cap\, \mathbb{Z}^{d}$ is the vertex set for the corresponding Delaunay polytope $D_{f} = \conv 
\sV(f)$. \ Conversely, if $D$ is a Delaunay polytope in $\mathbb{Z}^{d} \otimes \R$, there is a 
circumscribing empty ellipsoid  determined by a function $f_{D}\in \cP^{d}$. \ 
More precisely, since $D$ is Delaunay, there is a form $\varphi _{D}$, a center $\mathbf{c%
}$, and a radius $R$, so that $f_{D}(\mathbf{x})=\varphi_{D}[\mathbf{x}-%
\mathbf{c}_D]-R_D^{2}\le 0$ defines a circumscribing empty ellipsoid. \ Since  $f_{D}$ is 
positive on $\mathbb{Z}^{d}$, it is an element of $\cP^{d}$.

In the case where $\mathcal{E}_{f}$ is an empty cylinder, there is an infinite number of 
integer points lying on this surface. When this happens $\sV(f)$ is the set of all lattice 
points  on the boundary of a non-bounded \emph{polyhedron }$D_{f}=\conv \sV(f)$, which we  call 
a Delaunay polyhedron (we follow a recent convention that a polyhedron deserves to be called a   
a polytope when it is bounded). \emph{Although elements of $\sV(f)$ need not be vertices in the sense 
of real polyhedral geometry, we still refer to them as such and denote them by $\vertex D_f$.}   With this definition of
Delaunay polyhedron, for a given form 
 $\phi \ge 0$ Delaunay polyhedra need not form a tiling of $\mathbb{Z}^{d} \otimes \R$. They do 
form a face-to-face tiling if and only if the real kernel of the quadratic part of $f$ 
intersects with $\Z^d$ over a sublattice of the same dimension, i.e. if $\dim \ker_{\R} 
\phi=\dim  \{\ker_{\R} \phi \cap \Z^d\}$. Some authors refer to  this conditions as that $\phi$ 
has a rational kernel or as that $\phi$  has a rational radical, although both expressions are 
somewhat ambiguous -- a more accurate way of expressing this condition is to say that $\phi$ 
has fully rational kernel.
\begin{definition}
A polyhedron $D \subset \mathbb{Z}^{d} \otimes \R$ is called a Delaunay polyhedron, if there is  
$f \in \cP^{d} $ such that $D=\conv \sV(f)$. 
\end{definition}

\begin{definition}
Let $D$ be a Delaunay polyhedron in $\mathbb{Z}^{d} \otimes \R$. \ Then, the (inhomogeneous) 
domain $\cP_{D}^d$\ for $D$ is:
\begin{equation*}
\cP^{d}_{D}=\{\,f \in \cP^{d} \;\; \vline \;\;\text{ }D_{f}=D. \,\}
\end{equation*}
\end{definition}
\medskip

\noindent Such domains are relatively open convex cones that partition~$\cP^{d}$.

Each element $f\in \cP^{d}_{D}$ satisfies the linear equations $f(\mathbf{z})~=~0,\; \mathbf{z} 
\in D~\cap~\mathbb{Z}^{d}$. When $D$ is a single lattice point  and $d>0$, $\cP^{d}_{D}$  is a 
relatively open cone of dimension $\binom{d+2}{2}$. When $D$ is a 1-polytope, $\cP^{d}_{D}$ is 
a relatively open facet of the partition with dimension $\dim \cP^{d}-1=$ $\binom{d+2}{2}-1$. 
When the rank of the system $\{f(\mathbf{z})~=~0\; \vline \;\mathbf{z} \in 
D\cap\mathbb{Z}^{d}\}$ is full, i.e. $\binom{d+2}{2}$, the only possible $f$ is $0$ and 
$D=\mathbb{Z}^{d} \otimes \R$.  If $D$ is a $d$-dimensional ($d>0$) Delaunay simplex then $\dim \cP^{d}_{D}=$ $%
\binom{d+1}{2}$, but if $D$ is a perfect Delaunay polyhedron and $d>0$, then $\dim 
\cP^{d}_{D}=1$.

\begin{definition}
A function $p\in \cP^{d}$ is perfect if the system of equations 
\[ \left\{f(\mathbf{z})=\underset{\Z^d}{\min p}\;\; \vline\;\; \mathbf{z}\in \sV(p)\right\}\] has a unique 
solution $f=p$. \ In this case we also call $\conv \sV(p)$ a perfect polyhedron.
\end{definition}
The elements of 0- and 1-dimensional inhomogeneous domains are perfect, and the Delaunay 
polyhedra that determine such domains are perfect. The domain that consists of strictly 
positive constant functions corresponds to the empty Delaunay polytope. The perfect subsets 
$\sV(p)$ must be maximal among the subsets $\{\, \sV(f)  \;\; \vline \;\; f\in \cP^{d} \,\}$. 
With the exception of the 0-dimensional perfect domains, all inhomogeneous domains 
$\cP^{d}_{D}$ have proper faces that are inhomogeneous domains of lesser dimensions. Perfect 
quadratic functions are analogues of perfect  forms - both achieve their arithmetic minimum  at 
a sufficient number of points so that the representations of the minimum uniquely determine the 
polynomial.\

\   The following theorem and the next subsection show the important role played by perfect 
Delaunay polytopes in the theory of lattice Delaunay polytopes and Voronoi's L-types. 
\begin{theorem} (Erdahl, 1992)
Let $D$ be a Delaunay polyhedron in $\mathbb{Z}^{d} \otimes \R$.  Then
\begin{equation*}
\cP^{d}_{D}=\left\{\, \sum_{\{p\;  | \; \text{$p$ \emph{is perfect},}  \vert D \subseteq \sV(p) 
\}}\omega _{p}p \;\; \vline \;\; \omega_p \in \R_{>0} \,\right\},
\end{equation*}
\end{theorem}

This theorem shows that an arbitrary element $f\in \cP^{d}_{D}$ has the following 
representation:
\begin{equation*}
f=\sum_{\{p\;  | \; \text{$p$ is perfect,} \vert D \subseteq \sV(p) \}}\omega_{f,p}p,
\end{equation*}
where $\omega_{f,p} >0.$ 

\medskip

\noindent \textbf{The homogeneous domain of a Delaunay tiling: }\ Voronoi's classification 
theory for lattices, his \emph{theory of lattice types} (L-types), was formulated using 
positive definite forms and the fixed lattice $\mathbb{Z}^{d}$. \ In this theory two lattices 
are considered to be the same type if their Delaunay tilings are affinely equivalent. This 
characterization is not apparent in Voronoi's memoirs and was brought to light by Delaunay. 
Here is a brief reminder of Delaunay's approach to Voronoi's theory of L-types. Consider a 
positive definite quadratic form $\varphi $. \ Then a $\Z^d$-polytope $D$ is Delaunay relative 
to $\varphi $ if (1) it can be circumscribed by so called \emph{empty ellipsoid }$\mathcal{E}$ 
defined by an inequality  of the form
\begin{equation*}
\varphi [\mathbf{x-c}] \le R^{2}\text{,}
\end{equation*}
with center  $\mathbf{c}\in \mathbb{R}^{d}$ and radius $R \in \R_{>0}$,
 (2) $\mathcal{E}$ has no interior $\mathbb{Z}^{d}$%
-elements, and (3) $\vert D = \mathcal{E} \cap \mathbb{Z}^{d}$. \ The collection of all such 
Delaunay polytopes fit together facet-to-facet to tile $\mathbb{R}^{d}$, a tiling that is 
uniquely
determined by $\varphi $. \ This is the Delaunay tiling $\mathcal{D}%
_{\varphi }$ of $\mathbb{Z}^{d} \otimes \R$ relative to the form form $\varphi $.
If a second form $\vartheta $ has Delaunay tiling $\mathcal{D}%
_{\vartheta }$, and if $\mathcal{D}_{\vartheta }$ is
$GL(n,\Z)$-equivalent to $\mathcal{D}_{\varphi }$, then $\varphi $ and $%
\vartheta $ are forms of the same L-type. \

The description we give below requires that certain degenerate positive forms be admitted into 
the discussion, namely, those forms $\varphi $ for which $\ker_{\R} \varphi$ is a fully 
rational subspace of $\mathbb{Z}^{d} \otimes \R$.  The Delaunay polyhedra for such a form are 
themselves degenerate -- they are cylinders with axis $\ker_{\R} \varphi$ and Delaunay 
polytopes
as bases. \ These cylinders fit together to form the (degenerate) Delaunay tiling $%
\mathcal{D}_{\varphi }$.  For example, if $\mathbf{p} \in \mathbb{Z}^{d}$ is primitive, i.e. $\GCD(p_1,\ldots,p_d)=1$, 
then $\varphi [\mathbf{x}]=(\mathbf{p\cdot x})^{2}$ is such a form -- the kernel $\ker_{\R} \varphi$  
is the solution set for $\varphi [\mathbf{x}]=0$, and given by $\mathbf{p}^{\perp }$, which is 
fully rational. \ The Delaunay tiles are 
infinite slabs, each squeezed in between a pair of hyperplanes $\mathbf{p\cdot x=}k%
\mathbf{,p\cdot x=}k+1$, $k\in \mathbb{Z}$. \ These fit together to tile $\mathbb{Z}^{d} 
\otimes \R$. 

 Let $\Phi^d$ be the cone of positive definite forms and semidefinite forms with fully 
 rational kernels ($\Phi^d$ is sometimes referred to as the rational closure of the cone 
 of positive definite forms).  For each form 
$\varphi \in \Phi^d$ there is a Delaunay tiling ${Del}(\varphi)$ for $\mathbb{Z}^{d} \otimes 
\R$.

\begin{definition} If $\mathcal{D}$ is a Delaunay tiling for $\mathbb{Z}^{d}$, then
the following cone of positive definite quadratic forms
\begin{equation*}
\Phi _{\mathcal{D}}=\left\{\, \varphi \in \Phi^d \;\; \vline  \;\; Del(\varphi)=%
\mathcal{D}\, \right\}
\end{equation*}
is called an L-type domain.
\end{definition}

\noindent For this definition the Delaunay tilings can be the usual ones,
where the tiles are Delaunay polytopes -- or they could be degenerate
Delaunay tilings where the tiles are cylinders with a common axis $K$.

The relatively open faces of an L-type domain, are L-type domains. \ If $\mathcal{D}$ is a 
trangulation, $\Phi _{ \mathcal{D}}$ has full dimension $\binom{d+1}{2}$; this is the generic 
case. \ If $\mathcal{D}$ is not a triangulation, $\dim \Phi _{\mathcal{D}}$ is less than 
$\binom{ d+1}{2}$, and $\Phi _{\mathcal{D}}$ is a boundary cell of a full-dimensional L-type 
domain. \ In the case where $ \dim \Phi _{\mathcal{D}}=1$, the elements $\varphi \in \Phi 
_{\mathcal{D}}$ are called edge forms -- Voronoi showed that such domains are exactly extreme 
rays of full-dimensional L-type domains.

\par \emph{Let $\pi _{\Phi}$ be the projection operator from $\cP^{d}$ onto its quadratic part.} If 
$D \in \mathcal{D}$ and $\dim D=d$, then 
 $\Phi _{\mathcal{D}}\subset \pi_{\Phi}(\cP^{d}_{D})$. \ \ Since this containment holds for all 
Delaunay $d$-tiles $D\in \mathcal{D}$, there is the following description of $\Phi 
_{\mathcal{D}}$ in terms of inhomogeneous domains:
\begin{equation*}
\Phi _{\mathcal{D}}=\bigcap_{D\in \mathcal{D}}\pi _{\Phi}(\cP^{d}_{D}),
\end{equation*}
where the intersection is over all $d$-dimensional Delaunay polyhedra in $ \mathcal{D}$. (It is 
also true that the intersection can be taken over all Delaunay polyhedra in $ \mathcal{D}$.) \ 
Since the containment $\Phi_{\mathcal{D}}\subset \pi_{\Phi}(\cP^{d}_{D})$ also holds for all 
Delaunay tilings $\mathcal{D}$ that contain $ D$, there is the following description of 
$\pi_{\Phi}(\cP^{d}_{D})$ in terms of homogeneous domains:
\begin{equation*}
\pi_{\Phi}(\cP^{d}_{D})=\bigsqcup_{\mathcal{D\ni }D}\Phi _{\mathcal{D}},
\end{equation*}
where the disjoint union is over all Delaunay tilings of $\Z^d \otimes \R$ that contain $D$. \ 
This holds not only for full-dimensional  cells of $\cP^{d}_{D}$, but for cells of all 
dimensions. \ The last equality shows that $\pi _{\Phi}(\cP^{d}_{D})$ is tiled by L-type 
domains. It also establishes the following
\begin{proposition}
If $p\in \cP^{d}$ is perfect and, therefore, an edge form for an inhomogeneous domain, then 
$\pi _{\Phi }p$ is an edge form for an L-type domain. \ \
\end{proposition}
\noindent This result can also be established in a more direct way by appealing to the 
definition of L-type domain: if $D$ is a perfect Delaunay polyhedron, then $\pi _{\Phi 
}(\cP^{d}_{D})$ is a one-dimensional L-type domain. \ By definition, the elements of $\pi 
_{\Phi }(\cP^{d}_{D})$ are then edge forms. \

The converse of theorem does not hold - there are edge forms for L-type domains that are not 
inherited from perfect elements in $\cP^{d}$. \ \ Evidence has accumulated recently that the 
growth of numbers of types of edge forms with dimension is very rapid starting in six 
dimensions (see Dutour and Vallentin, 2003), but the growth of perfect inhomogeneous forms
is much less rapid - there is some hope that a complete classification can be made through 
dimension nine.
 
\section{\large Inhomogeneous domains and L-types\\ for arbitrary lattices}

In the previous section we only treated the case of  ``ground lattice" $\Z^d$ and tilings of 
 $\Z^d \otimes \R$. Often, it is more convenient to work with a ground lattice which is a 
 linear or affine section of $\Z^d$ or even  a centering (superlattice) of $\Z^d$.   A point lattice in a real or rational  affine space 
is a set $\LL$ such that for any $x \in \LL$ the set of vectors $\{y-x\;\vline\; y \in \LL\}$ is a free $\Z$-module.  $\Z^d$ is a $\Z$-module, i.e. a ``vector" lattice, but the ``ends of vectors" of $\Z^d$ form a point lattice in $\Z^d \otimes \R$ -- it is customary to use the symbol $\Z^d$ in both senses. If $\Lambda$ is a section of $\Z^d$ defined by an affine equation, where the $\Z^d$  is considered as a $\Z$-module, then the set  $\Lambda-\Lambda=\{\z-\z'\;\;\vline\;\;\z,\z' \in \Lambda\}$ is a module; if $\LL$ is an affine section of point lattice $\Z^d$, then $\LL$ is a point lattice. The point lattice setup is more intuitive in the study of quadratic \emph{functions} and Delaunay polytopes, while the ``vector" ($\Z$-module) setup is more natural for studying quadratic \emph{forms} and, perhaps,  Dirichlet-Voronoi  cells. A case can be made for either setup for Delaunay tilings, for they show both homogeneous and inhomogeneous aspects. We will use the term "affine sublattice" in both situations; the distinction will be clear from the context or notation. 

For example, one convenient way to describe the root
 lattice $A_n$ is to start from $\Z^{n+1}$ and consider a sublattice of $\Z^{n+1}$ 
 defined by equation $z_1+\dots+z_{n+1}=1$ with respect to the Euclidean lattice norm 
 $\parallel~\parallel$
 (sum of squares) on $\Z^{n+1}$.
Denote by $\Lambda$ the lattice of all vectors $\z-\z'$ where both $\z,\z'$  lie in $\LL=\{\x \in 
\Z^{n+1}\;\vline\;\sum_1^{n+1} x_i=1\}$.
 The resulting pair $[\Lambda,\parallel~\parallel\:]$ is 
 known as the root lattice $A_n$. It can also be regarded as a point lattice $\left(\LL, \underset{1}{\overset{n+1}{\sum}} x_i^2 |_{\LL}\right)$.
 
 A pair $[\Lambda,\phi]$, where $\Lambda$ is a free $\Z$-module of finite rank and 
 $\phi$ is a quadratic form on   $\Lambda$ is called a \emph{quadratic lattice} ( \emph{quadratic module} is also used).  A pair $(\LL,f)$, where $\LL$ is a point lattice of finite rank and $f$ is a quadratic function on $\LL$ is called a\emph{ quadratic point lattice. } Quadratic lattices 
 $[\Lambda,\phi]$ and $[\Lambda_1,\phi_1]$ are called 
isomorphic if there is a $\Z$-module isomorphism $m:\Lambda \rightarrow \Lambda_1$ such that 
$\phi[\x]=\phi_1[m\x]$ for any $\x \in \Lambda$.  Quadratic point lattices 
 $(\LL,f)$ and $(\LL_1,f_1)$ are called 
isomorphic if there is a point lattice isomorphism $m:\LL \rightarrow \LL_1$ such that 
$\pi_{\Phi}f[\x]=\pi_{\Phi}f_1[m\x]$ for any $\x \in \LL$.  
Although $\LL$ does not have a preferred origin, the quadratic part of  $f|_{\LL}$ is independent of the choice of such.
The notion of Delanay tiling makes sense for both types of lattices, the only difference is that in the case of quadratic "vector" lattices it has a preferred origin. In the former case we will denote 
the Delaunay tiling of $\Lambda$ with respect to $\phi$  by $Del(\Lambda,\phi)$, and in the latter case we denote the Delaunay tiling of $\LL$ with respect to the quadratic part of $f$ by $Del(\LL,\pi_{\Phi}f)$; when $\phi$ or $f$ are defined on a superset of a lattice we may omit the explicit respriction in the notation, e.g., we may use $\phi$ instead of more proper $\pi_{\Phi}\phi|_{\LL}$.

\section{\large Structure of Perfect Delaunay Polyhedra}
In this section we will describe how to construct new arithmetic types of perfect Delaunay
polytopes in dimension $d$ from known  types of perfect Delaunay polytopes in dimension $d-1$ 
for $d>2$. Since there are no perfect Delaunay $d$-polytopes for $1<d<6$, this construction can 
be used starting from $d=7$.

 Below $\phi$ denotes any positive quadratic 
form, i.e. $\phi[\z] \ge 0$ for any $\z$.  Recall that we denote by $\pi_{\Phi} p$ the 
quadratic form part of a polynomial $p$. Our starting point is the following structural 
characterization. 

\begin{theorem}(Erdahl, 1992) A polyhedron $P \in Del(\Lambda,\phi)$ is perfect
 if and only if there is a quadratic \emph{function }$p$ on $\Lambda$ with $\pi_{\Phi} p=\phi$ such 
 that
 $P \cap \Lambda=\sV(p)=\{\vv + \z \; \vline \; \vv \in \vertex D,\: \z \in \Gamma\}$, where 
$D$ is a perfect \emph{polytope } in $Del(\Lambda \cap \aff D,\pi_{\Phi}(\phi|_{\aff D}))$ and $\Gamma$ is 
a submodule of $\Lambda$ such that $\Lambda$ is the direct sum of modules $(\Lambda \cap \aff 
D)-(\Lambda \cap \aff D)$ and  $\Gamma$ (or, equivalently, point lattice $\Lambda$ is the direct affine sum of point set  $\Lambda \cap \aff 
D$ and   module $\Gamma$). 
\end{theorem}

As an example, consider primitive vectors $ \mathbf{a,b} \in \Lambda=\mathbb{Z}^{d}$ \ such 
that $\mathbf{a\cdot b=}1$. Let $\phi[\z]=(\mathbf{a}\cdot \z)^2$ and let us use 
$\mathbf{a}^{\perp }$ for $ \left\{ \x \in \R^d \: | \: \x \cdot \mathbf{a} = 0 \right\}$. Then 
$\mathcal{E}_{0}=\conv\{\mathbf{0,b}\}$ is an empty perfect ellipsoid in the 1-lattice $(\aff%
\mathcal{E}_{0})\cap \mathbb{Z}^{d}$. The set $\{ \mathbf{0,b} \}= \mathcal{E}_{0} \cap 
\mathbb{Z}^{d}$ is the vertex set for a perfect Delaunay \emph{polytope} $D$ in $\aff 
(\mathcal{E}_{0} \cap \mathbb{Z}^{d})$, and $\Gamma=(\mathbf{a}^{\perp } \cap \mathbb{Z}^{d})$. 
We have $\mathbb{Z}^{d}=\Z\b \oplus (\mathbf{a}^{\perp } \cap \mathbb{Z}^{d})$.  Thus, 
 the polyhedron $P$  is Euclidean direct product of $\conv\{\mathbf{0,b}\}$ and $ \mathbf{a}^{\perp }$. The perfect 
quadric circumscribed about $P$ is given, up to a scaling factor, by the equation $p(\x)=0$, 
where  $p(\mathbf{x})=(\mathbf{a\cdot x})( \mathbf{a \cdot x-}1)=\phi[\x]- \mathbf{a \cdot x}$. 
\ The preimages of negative real values of $p$ lie between two hyperplanes $\mathbf{a\cdot 
x=}0$ and $\mathbf{a\cdot x=}1 $, a region whose closure is a degenerate perfect Delaunay 
polyhedron $P$. 

By the above theorem, for each perfect $P \in Del(\Lambda,\phi)$  the set $\vert P$ of lattice points on $P$ can be written as $ \vertex D \oplus \Gamma$, where $D$ is a Delaunay polytope and $\Gamma$ is a (vector) sublattice of $\Lambda$. The direct sum $\oplus$ sign means that (i) each 
vertex of $P$ can be represented as a vertex of $D$ plus a vector of $\Gamma$ and (ii) \emph{for any fixed choice of the pair $(D, \Gamma)$} such 
representation is unique.  We call this a direct decomposition of $P$. If $P$ is a polytope, it 
is obviously unique. What are the possibilities in the case of a degenerate $P$?

\begin{lemma} 
Suppose $\vertex P={\vertex D} \oplus \Gamma=\vertex D' \oplus \Gamma$ are direct 
decompositions of a perfect polyhedron $P$ in $Del(\Lambda,\phi)$, where $D$ and $D'$ are Delaunay polytopes and $\Gamma$ is a (vector) sublattice of $\Lambda$.  If $\aff D $ is parallel to $\aff D'$, then $D'$ is a translate of $D$ by a vector of $\Gamma$.
\end{lemma}
\begin{proof} Since
 $\vertex P=\vertex D \oplus \Gamma$, any $v' \in \vertex D'$ can be written in a unique way as 
$v+\uu$, with $v \in \vertex D$ and $\uu \in \Gamma$. 
Since $\aff D'$ is parallel to $\aff D$,  we have $\aff D'= \aff D + \uu$ and 
$D+\uu \subset \aff D'$. In particular, $ \vertex D+\uu \subset \aff D' \cap \Lambda$. So, 
 $\vertex D+\uu \subset \vertex P \cap \aff D'$. 
Since 
 $\vertex P=\vertex D' \oplus \Gamma$,   
 we have $\vertex P \cap \aff D'=\vertex D'+0=\vertex D'$ and $\vertex D+\uu \subset \vertex D'$.
   Since $D$ is a perfect Delaunay polytope in 
 $Del(\aff D \cap \Lambda,\pi_{\Phi}(\phi|_{\aff D}))$, its translate $D+\uu$ is also a perfect Delaunay 
 polytope 
 in $Del(\aff D' \cap \Lambda,\pi_{\Phi}(\phi|_{\aff 
 D'}))$. But this means $D+\uu$ is maximal in $Del(\aff D' \cap \Lambda,\pi_{\Phi}(\phi|_{\aff 
 D'}))$ and $D'=D+\uu$.
\end{proof}

The following theorem gives a universal construction of a perfect Delaunay polytope in 
dimension $d+1$ from a perfect Delaunay polytope in dimension $d$.

\begin{theorem}\label{thm:step}Let $P$ be a perfect \emph{polytope} in $Del(\Lambda,\phi)$ and let 
$p$ be its \emph{perfect} quadratic function, i.e. $\vertex P=\sV(p)$ and $\phi=\pi_{\Phi} p$. Suppose 
$D \in Del(\Lambda,\phi)$ is another Delaunay cell of full dimension, which is not a 
$\Lambda$-translate of $P$. If $\e \notin \Lambda$, then there is a positive definite form 
$\psi$ on $\Lambda \oplus \Z \e$ and a perfect polytope  $P'$ in $Del(\Lambda \oplus \Z \e, 
\psi)$ such that $P'\cap \aff \Lambda=P$ and $P' \cap \{\aff \Lambda+\e\}=D+\e$
\end{theorem} 
\begin{proof}
We can safely assume $\mathbf{0} \in \vertex P$ and $\mathbf{0} \in \vertex D$. Let then $f(\y)=\psi[\y] +l(\y)$ be   a quadratic function with linear part $l$ such that:  \texttt{(i)} $f(\y)=0$ is the equation of a quadric passing through the vertices of $P$ and $D + \e$, \texttt{ (ii)} $f|_{\Lambda}=p$. 
Since 
any $\y \in \Lambda \oplus \Z \e$ can be uniquely written as $\x + k\e$ with $\x \in \Lambda$ and $k 
\in \Z$, we have $\psi[\y]=\psi[\x]+\psi[k\e]+2\psi(\x,k\e)= 
\phi[\x]+k^2\psi[\e]+2k\psi(\x,\e)$ and $l(\y)=l(\x)+kl(\e)$. $l|_{\Lambda}$ is known from $p$. So, $f$ is determined uniquely if and only if  $\psi[\e]$, $\psi(\e,\e_1),\ldots,\psi(\e,\e_d)$, $l(\e)$ are fixed. 

We will first show that for any positive ($\ge 0$) value of $\psi[\e]$, function $f$, that satisfies the conditions \texttt{(i)}, \texttt{(ii)}  is unique, and $f(\y) \le 0$ defines a quadric, circumscribed about $P \cup \{D+\e\}$, such that no points of $\Lambda \oplus \Z \e$ satisfy $f(\y) <0$ . Let us fix some $\psi[\e]\ge0$. Since $\dim D=d$, it has $d+1$ affinely independent vertices. For each non-zero $\vv \in \vertex D$ consider equation $f(\vv+\e)=0$. If $\vv_1,\ldots,\vv_d$ are independent, then the resulting system of equations is independent. (The origin has already been used to set the constant term of $f$ to $0$.) We know the line joining the circumcenters  of $D+\e$ and $P$ is perpendicular to $\Lambda$ with respect to $\psi$. This condition can be written as $\psi(\cc_D-(\cc_P+\e),\vv)=0$, where $\vv \in \vertex D \diagdown \mathbf{0}$, and $\cc_D$, $\cc_P$ are the circumcenters of $D$ and $P$. The resulting system of $d+1$ linear equations has full rank.  It is inhomogeneous because $\phi[\vv]+\psi[\e]>0$. So, the  solution is unique. Positivity of $\psi$ follows from positivity of $\phi$, $\psi[\e]\ge0$ and elementary geometric considerations. For a sufficiently small value of $\psi[\e]\ge0$ the interior of   $\{\z \in  (\Lambda \oplus \Z \e) \otimes \R \; \vline\;f(\z)\le 0\}$ is free of lattice points. This is easy to establish by a calculation, or just by observing that large values of $\psi[\e]$ correspond to ellpsoids that lie inside of an infinite slab of $(\Lambda \oplus \Z \e) \otimes \R$, which is squeezed between hyperplanes  corresponding to $k=-1$ and $k=2$. 

Thus, under the assumptions of  the theorem, for a fixed $\psi[\e] \ge 0$,  the quadric  circumscribed about  $P \cup \{D+\e\}$ is  uniquely defined, has positive quadratic part, and is Delaunay for  sufficiently large $\psi[\e]$'s.  We can now start continuously decreasing $\psi[\e]$, at a constant rate, until the ellipsoid hits a new lattice point, i.e. a point of $\Lambda \oplus \Z \e~\diagdown~(P \cup \{D+\e\})$. If it does hit a new lattice element $\a$ in a finite time, we are done, since $f(\a)=0$ will give us a linear  equation independent of those that we used to determine    $\psi(\e,\e_1),\ldots,\psi(\e,\e_d)$ and $l(\e)$. This new equation will determine $\psi[\e]$; the resulting Delaunay ellipsoid $f(\z)=0$ will be perfect.  Note that  $f(\z)=0$ may contain other new lattice elements  besides $\a$. If such collision does not happen in finite time, then $\psi[\e]=0$ determines a perfect Delaunay surface, which is a cylinder with elliptic base. By the preceeding lemma any section of this cylinder by a sublattice parallel to $\Lambda$ is the vertex set of a lattice translate of $P$ (by a vector $n\e$ for some $n \in \Z$). But the section of  $\Lambda \oplus \Z \e$ given by $k=1$ is a translate of $D$, which is not a translate of $P$. It only means that the infinite cylinder scenario is impossible and there exists a required $f$.
\end{proof}
 
 Since $P$ is perfect in  $Del(\Lambda,\phi)$, it determines $\phi$ (up to a scale factor) and  the corresponding Delaunay tiling of $\Lambda$ uniquely. However, there can be different choices of $D$. For example, the 35-tope of Erdahl and Rybnikov lives in a 7-dimensional Delaunay tiling which has many  inequivalent (even combinatorially!) types of Delaunay polytopes of full dimension. Applying this theorem with different choices of $D$ results, in general but not always, in arithmetically distinct perfect polytopes in the next dimension. For example, the Delaunay tiling determined by $P=G_7$ (this is the Delaunay tiling of $[E_7,\parallel~\parallel]$) consists of translates of a single copy of $G_7$, and  many (more than two -- one always has at least 2 for a simplex) translation classes  of Delaunay simplices of index 2 (see Erdahl 1992 for an exact description). All these simplices are arithmetically equivalent, in particular, they are isometric relative to the form determined by $G_7$; however, there is more than one orbit of such simplicies with respect to the subgroup generated by lattice translations and inversion with respect to the origin. But no matter which copy of such a simplex we use for $D$, the resulting perfect Delaunay  polytope and quadratic lattice in dimension 8 is the same. On the other hand, different choices of translationally-inequivalent $D$'s in the case of the Delaunay tiling determined by the 35-tope sometimes lead to  inequivalent perfect Delaunay polytopes in dimension 8.

So, this construction always works, but is not always uniquely determined by $P$. The only known cases where one cannot find a $D$ which is not a translate of $P$, are those of $d=0$ and $d=1$. A surprising fact is that for $d=0$ the construction works even without such a $D$, since $[0,1]$ is perfect Delaunay in $\Z^{d+1}$ for $d=0$. We do not know of any other example of a perfect Delaunay polytope which can be translated  by a lattice vector to every Delaunay polytope of full dimension in the Delaunay tiling! Very likely that such an example does not exist. If we drop the perfection condition, the only example that we know of is a $d$-parallelepiped. There should be  a way to  show that no other examples exist, but, embarassingly,  we do not have a proof at the moment.
 
\section{\large Infinite Sequences of \\ Perfect Delaunay polytopes}
To simplify notation we will denote an affine sublattice spanned by all affine integral combinations of a set of lattice points $S$ by $\aff_{\Z}S$, and we will use the same abbreviation for the (vector) lattice spanned by all integral linear combinations of vectors $\{\z-\z' \; \vline\;\z,\z' \in S\}$.   
Denote by $\mathbf{j}$ the vector of all ones, i.e. $[1,\ldots,1]$, in $\Z^d$. Denote by $\mathbf{V}_{s,2k}^{d}$  the following set in $\mathbb{Q}^{d}$:
\begin{equation*}
\mathbf{V}_{s,2k}^{d}\triangleq\{\, [1^{s},0^{d-s}]-\frac{s-1}{d-2k} \mathbf{j}\,\} \times \binom{d}{s},
\end{equation*}
where $s,k \in \N=\{0,1,\ldots \}$.  All permutations of entries are taken, so that $|\mathbf{D}_{s,2k}^{d}|=\binom{d}{s
}+\binom{d}{s+1}=\binom{d+1}{s}$. In general, whenever we use a shorthand notation for vectors or points such as $[1^{s},0^{d-s}]$, \emph{it is understood that all permutations of the components are taken, except for those components that are separated from others by semicolumns on both sides,or by a semicolumn and a bracket}.
Also, let us set
\begin{equation*}
\mathbf{D}_{s,2k}^{d}\triangleq \{\,\vv-\vv'\;\;\vline\;\;\vv,\vv' \in \mathbf{V}_{s,2k}^d\,\} 
\end{equation*}
and 

\begin{equation*}
P_{s,2k}^{d}\triangleq \conv\left\{\;\mathbf{V}_{s,2k}^{d}\cup -\mathbf{V}_{s,2k}^{d}\cup \mathbf{V}_{s+1,2k}^{d}\cup -\mathbf{V}_{s+1,2k}^{d} \;\right\}.
\end{equation*}

The following is the main constructive theorem of
this paper.

\begin{theorem}\label{main_theorem}
Let $d,s,k \in \N$. If $2\le k$ and $1 \le s \le \frac{d}{2k}$, then $P_{s,2k}^{d}$
is a symmetric perfect Delaunay polytope for the quadratic lattice $[\Lambda _{s,2k}^{d}, \varphi _{s,2k}^{d}]$, where $\Lambda _{s,2k}^{d}=\aff_{\Z} \vertex P_{s,2k}^{d} $. The
circumscribing empty ellipsoid is defined as $\{\x \in \R^d\;\vline\;\varphi _{s,2k}^{d}[
\mathbf{x}]~\le~R^{2}\}$, where
\begin{equation*}\label{formula_for_phi}
\varphi_{s,2k}^{d}[\mathbf{x}]=4k(d-2sk-k)|\mathbf{x}|^{2}+\left(d^{2}-(4k+2s+1)d+4k(2s+k)\right)(\mathbf{j\cdot
x)}^{2}.
\end{equation*}
\end{theorem}
The origin is the center of symmetry for $P_{s,k}^{d}$ and does not belong to $\Lambda _{s,k}^{d}$. Thus, each pair of integers $(s,k)$, for $s\geq 1,k\geq 2$ and $s \le \frac{d}{2k}$, determines an infinite sequence of symmetric perfect Delaunay polytopes, one in each dimension. For $s=1,k=2$ the infinite sequence is the one
described in the opening commentary, i.e.,   $G_{7}^d, d\geq 7$, where the initial term is the 
Gosset polytope $G_7^7=G_7=3_{21}$ with $56$ vertices. \

The $\binom{d+1}{s}$ diagonal vectors $\mathbf{D}_{s,2k}^{d} \cup \mathbf{D}_{s+1,2k}^{d} $ for
$P_{s,2k}^{d}$ have the origin as a common mid-point, forming a
segment arrangement that generalizes the cross formed by the
diagonals of a cross-polytope. \ Moreover, these \ $\binom{d+1}{s}$
diagonals are primitive and belong to the same \emph{parity class}
for $\Lambda_{s,2k}^{d}$, namely, they are equivalent modulo
$2\Lambda_{s,2k}^{d}$. \
 \ \ More generally, primitive
vectors $ \uu ,\vv $ in some lattice $\Lambda $, with mid-points equivalent modulo 
$\Lambda$, are necessarily equivalent modulo $2\Lambda $ and thus belong to the same parity 
class. \ \ And conversely, the mid-points of lattice vectors $\uu ,\vv \in \Lambda $ belonging 
to the same parity class are equivalent modulo $\Lambda $. By analogy with the case of 
cross-polytopes, we call any such arrangement of segments or vectors a cross. The convex hulls 
of such crosses often appear as cells in Delaunay tilings -- cross polytopes are examples, as 
are the more spectacular symmetric perfect Delaunay polytopes.  \ There is a criterion, 
essentially  due to Voronoi (1908), but first formulated and formally proved by Baranovskii 
(1991), that determines whether a cross is Delaunay:  \emph{Let $\Lambda $ be a lattice, 
$\varphi $ a positive definite form, $C$  the convex hull of a cross of primitive vectors 
belonging to the same parity class.  Then $C$ is Delaunay relative to $\varphi $ if and only if 
the set of vectors forming the cross is the complete set of vectors of minimal length, relative 
to $\varphi$, in their parity class. } \ \ This is the criterion we have used to establish the 
Delaunay property for the symmetric perfect Delaunay polytopes $P_{s,2k}^{d}.$

The following result shows that asymmetric perfect Delaunay polytopes can
appear as sections of symmetric ones.

\begin{theorem}\label{sym_from_antisym}
For $d\geq 6$ let $\mathbf{u}=[-1^{2};1^{d-1}] \in \Z^{d+1}$. \ Then
\begin{equation*}
G_{6}^d=\conv\{\, \mathbf{v} \in \vert P_{1,4}^{d+1} \;\; |  \;\; \mathbf{v \cdot u}=\frac{1}{2} 
\, \}
\end{equation*}
is an asymmetric perfect Delaunay polytope for $(\Lambda _{1,4}^{d}, \psi _{1,4}^{d+1})$, where $\Mu_{1,4}^{d}=\aff_{\Z} \vertex G^d_d$. 

The
circumscribing empty ellipsoid is defined as 
$\{\x \in \aff_{\R}G_6^d\;\;\vline\;\; \psi _{1,4}^{d+1}[\mathbf{x}]~\le~R^{2}\}$, where
\begin{equation*}
\psi _{1,4}^{d+1}[\mathbf{x}]=8(d-5)|\mathbf{x}|^{2}+(d^{2}-9d+22)(\mathbf{%
j\cdot x)}^{2}.
\end{equation*}
\end{theorem}

\noindent This is the infinite sequence of $G_6$-topes described in the introduction, with  
Gossett polytope $G_6^6=G_6=2_{21}$ as the initial term for $d=6$. Note that the formula for the form on $\Mu _{1,4}^{d}$  is given with respect to  $\Z^{d+1}$.

The terms in this sequence have similar combinatorial properties. \ For example, the lattice 
vectors running between vertices all lie on the boundary, in all cases. These lattice vectors 
are either edges of simplicial facets, or diagonals of cross polytope facets--there are two 
types of facets, simplexes and cross polytopes. \ The Gossett polytope $G_6$ has $27$ 
5-dimensional cross-polytopal facets, but for other members of the series  $G_6^d$, the number of cross-polytopal facets is equal to 
$2d.$ $G_6$ can be found as a section of $G_6^7$, but $G^8_6$ and $G^9_6$ do not have sections 
arithmetically equivalent to $G_6$.

We summarize the properties of both  $G$-series in 
Table 1. Table 2  (note that entries separated by semicolumns are fixed) gives coordinates of the vertices of $G_6^d$-topes (named $\Upsilon^d$-polytopes 
there), discovered in 2001 by Erdahl and Rybnikov.

\begin{table}
\setlength{\extrarowheight}{4pt}
\[ \begin{array}{|c|c|c|c|}\hline
Polytope & \dim P & | \vert P| & Symmetry \\ \hline P_{1,4}^{d+1} & d+1 &
2\left({\binom{d+1}{1}}+{\binom{d+2}{1+1}}\right)=2{\binom{d+2}{2}} &
\textrm{centrally-symmetric} \\ \hline \Upsilon^{d}=G^{d}_6 & d &
\frac{d(d+2)}{2} - 1 & asymmetric \\ \hline
\end{array} \]
\caption{Properties of constructed perfect Delaunay polytopes}
\label{degree-4-number-table}
\end{table}

\begin{table}
\setlength{\extrarowheight}{5pt} 
\begin{center}
\begin{tabular}{|c|c|c|}\hline
$\mathbf{[ 0^{d}]}\times 1$ & $\mathbf{[-1,0^{d-2};1]}\times (d-1)$ & $
\mathbf{[1^{d-1};-(d-3)]}\times 1$ \\ \hline
$\mathbf{[0,1^{d-2};-(d-4)]}\times (d-1)$ & $\mathbf{[1^{2},0^{d-3};-1]}\times 
\frac{(d-1)(d-2)}{2}$ & $\mathbf{[1,0^{d-2};0]}\times (d-1)$ \\ \hline
\end{tabular}
\caption{Vertices of $G_6^d$-topes (same as $\Upsilon^d$) generalizing Gosset's $G_6$ on 27 vertices}
\end{center}
\end{table}

\subsection{\large Delaunay Property Part of Main Theorem} Let $\e_1,\dots,\e_d$ stand for the canonical basis of
${\mathbb Z}^d$. We will use the
following notation: $\phi_1[\x] = \left(\sum_{i=1}^d x_i\right)^2$, $\phi_2[\x] = \left| \x - 
\frac{\sum_{i=1}^d x_i}{d} {\bf j} \right|^2$; $\phi_2[\x]$ is the squared Euclidean distance 
from $\x$ to the line $\R {\bf j}$. The following proves the Delaunay property of 
polytopes $P_{s,2k}^d$ asserted by Theorem \ref{main_theorem}.

\begin{theorem}
\label{delaunay-property-theorem}
Let $d,s,k \in \N$, and let $2 \le k$ and $1 \le s \le \frac{d}{2k}$.
Let also \[L^1 = \{[{(-1)}^{k}, 1^{d-{k}}]\} \times \binom{d}{k}\]
\[\Lambda =
\Z\left\langle \e_1,\dots,\e_d,\frac{\bf j}{d-2k}\right\rangle\]
 \[\Lambda_1 =
\lbrace \z \in \Lambda\;\; \vline \;\; \z \cdot L^1 \equiv 1 \mod 2
\rbrace\]
 Then there is a positive definite quadratic form of the type
\begin{equation}\label{quadratic_form_equation} \phi[\x] = \alpha \phi_1[\x] + \beta
\phi_2[\x],
\end{equation}
with $\alpha, \beta \in \Q_{>0}$, such that $P_{s,2k}^d$ is a Delaunay
polytope in the lattice $(\Lambda_1, \phi|_{\Lambda_1})$.
\end{theorem}

Below $n=d-2k$

\begin{lemma} \label{structural_lemma}
Suppose $\phi[\x] = \alpha \phi_1[\x] + \beta \phi_2[\x]$ where $\alpha,
\beta > 0$. Then all points $\z\in\Lambda_1$ which are closest to $\mathbf{0}$
with respect to $\phi$, ie.,
$$ \phi[\z] = \min\lbrace \phi(u) : u\in\Lambda^0 \rbrace $$
are, up to permutations of components, of the type \[\sgn(l)[ 1^{|l|}; 0^{d-|l|} ] +
\frac{a}{n}{\bf j},\] with $-\frac{d}{2}  \le  l < \frac{d}{2}$ , for some
$a\in{\mathbb Z}$. Furthermore, each point of   $\Lambda_1$, closest to the origin, has only one representation of the described type.
\end{lemma}

\begin{proof} Suppose $\z \in \Lambda_1$ is minimal with respect to $\phi$, i.e., $\phi[\z]=\underset{\Lambda_1}{\min} \phi$. Let $\z = a_1\e_1 + \dots + a_d \e_d + a
\frac{\bf  j}{n}$, where $a_1,\dots,a_d,a \in {\mathbb Z}$ and $A = a_1 + \dots +
a_d$. We have
\begin{equation}
%\begin{split}
\phi_2[\z] = \phi_2\left[ \sum_{i = 1}^{d} a_i \e_i \right] =
 \left| \sum_{i = 1}^{d} a_i \e_i - A\frac{\mathbf  j}{d}\right|^2 = \sum_{i=1}^{d} a_i^2 - \frac{A^2}{d}.
%\end{split}
\end{equation}

Let us prove that the coefficients  $\{a_i\}$ are of two consecutive integer values. Suppose, to the contrary,  there
are $a_i$ and $a_j$ such that $a_i - a_j \ge 2$. Consider the vector
$\z' = a_1 \e_1 + \dots + (a_i - 1)\e_i + \dots + (a_j + 1) \e_j +
\dots + a_d \e_d + a\frac{\mathbf  j}{n}$. We have $\phi_2[\z'] -
\phi_2[\z] = (a_i - 1)^2 + (a_j + 1)^2 - a_i^2 - a_j^2 = 2(a_j - a_i)
+ 2 \le -2$ and ${\mathbf  j}\cdot \z' = {\mathbf  j} \cdot \z$. Since
$\z' \in \Lambda_1$ and $\phi[\z'] < \phi[\z]$, it follows
that the vector $\z$ is not closest to $\mathbf{\z}$ which is a contradiction.

Now, let $b$ be the smallest of the values of the coefficients $\{a_i\}$. Subtract
$b\mathbf{j}$ from the first part and add an equal value of
$bk\frac{\bf j}{n}$ to the second part of the existing representation of
$\z$ as an integral linear combination of $\e_i$'s and $\frac{\bf j}{n}$. After a permutation of the components $\z$ is
equal to $\sgn(l)[1^{|l|},0^{d-|l|}] + (a + bk)\frac{\bf j}{n}$ where $0\le l < d$. If
$l \ge \frac{d}{2}$, again subtract $[1^d]$ from the first summand and add
${\bf j}$ to the second summand to get the required representation.

Note that although $\dim \Lambda=d$, we represent vectors of $\Lambda_1$ as integral linear combinations of $d+1$ vectors; thus, we have to understand if the representation described in the statement of the lemma is unique for each of the points of $\Lambda_1$ that are closest to the origin. To prove that each minimal vector of $\Lambda^0$, with respect to $\phi$, has only one encoding of the described type, note that the components  of the
vector $\sgn(l)[ 1^{|l|}; 0^{d-|l|} ] + a\frac{{\bf j}}{n}$ are of at most two
values. In $\sgn(l)[ 1^{|l|}; 0^{d-|l|}]$, either $1$'s fill the
positions of the  largest value, or $-1$'s fill the positions
of the smallest value. Since $-\frac{d}{2} \le l <
\frac{d}{2}$, the choice is unique. 
\end{proof}

\begin{lemma}
\label{2nd_structural_lemma} All points of the affine lattice
$\Lambda_1$, which are closest to $\mathbf{0}$ with respect to a form $\phi =
\alpha \phi_1 + \beta \phi_2$, where $\alpha,\beta > 0$, belong to
the set $\lbrace \z \in \Lambda_1 \; \;\vline\; \; \left| \z
\cdot {\bf j} \right| \le \frac{d}{n}\rbrace$. If $\z$ is
closest to $\mathbf{0}$ with respect to $\phi$ and  $\left| \z \cdot
{\bf j} \right| = \frac{d}{n}$, then $\z\in\lbrace\pm\frac{\bf
j}{n}\rbrace$.
\end{lemma}

\begin{proof} By multiplying both $\alpha$ and $\beta$ by the same positive number, we can assume that
the minimal value of $\phi$ on $\Lambda_1$ is $1$. Since $\frac{\bf j}{n}
\in \Lambda_1$, $\phi[\frac{\bf j}{n}] = \alpha (\frac{d}{n})^2 \ge 1$.
Let $\z\in{\mathbb R}^d$ be a point with $\phi[\z]=1$. Represent
$\z = \frac{\gamma}{n} {\bf j} + \uu$ where $\gamma\in{\mathbb R}$,
$\uu\cdot{\bf j} = 0$. We have $\phi[\z] = \alpha\frac{\gamma^2
d^2}{n^2} + \beta |\uu|^2 = 1$, therefore $\alpha\frac{\gamma^2 d^2}{n^2}
\le 1$. Since $\alpha(\frac{d}{n})^2 \ge 1$, we have proven that $|\gamma|
\le 1$, so $|\z\cdot{\bf j}| = |\gamma|\frac{d}{n} \le \frac{d}{n}$.

If the latter inequality holds strictly, then $|\gamma|=1$ and
$\phi[\z] = \alpha\frac{d^2}{n^2} + \beta |\uu|^2 = 1$. Since $\alpha
(\frac{d}{n})^2 \ge 1$, we necessarily have $\beta = 0$ so $\z =
\pm\frac{\bf j}{n}$.
\end{proof}

\medskip

\par \noindent \textbf{Proof of Theorem 6:}
Let  $\phi_{\alpha,\beta}$ be the type (\ref{quadratic_form_equation}), i.e. $\phi_{\alpha,\beta}[\x]= \alpha \phi_1[\x] + \beta \phi_2[\x]$.  Let us first describe a ``reduced" set of points $M$, such that any point of
$\Lambda_1$  where the minimum of   $\phi_{\alpha,\beta}$ over $\Lambda_1$ is attained,  belongs to this set, possibly after a permutation of the components. Then we "handpick" the form $\phi=\phi_{\alpha,\beta}$ so that the vertices of $P^d_{s,2k}$ are the minimal
\emph{points} and other points of $\Lambda_1$ are not. Below is an implementation of
this plan.

\medskip

Consider the set
\begin{equation}
M = \left\lbrace\z \in  \Lambda_1  \; \;\vline \;\; \z=\sgn(l)[ 1^{|l|} ;0^{d-|l|} ] + a\frac{\bf j}{n},   0 \le \z\cdot{\bf j} < \frac{d}{n},\; -\frac{d}{2} \le l < \frac{d}{2} \right\rbrace \cup \lbrace \frac{\bf j}{n} \rbrace.
\end{equation}
By the two previous lemmas,  every  $\z \in\Lambda_1$, which is closest to the origin of $\Lambda$ with respect to $\phi$, perhaps after a
permutation the components, belongs to this set.

For $\z =\sgn(l) [ 1^{|l|}; 0^{d-|l|} ] + a\frac{\bf j}{n}$, we have
\begin{equation}
\begin{split}
\z\cdot{\bf j} = l + \frac{ad}{n}, \\
\phi_1[\z] = \left(l + \frac{ad}{n}\right)^2, \\
\phi_2[\z] = |l| - \frac{l^2}{d},\\
L^1\cdot\z\equiv a+l \!\!\mod 2.
\end{split}
\end{equation}

From the calculations above, we get that $M$
\begin{equation}
\label{M_equation} M = \left\lbrace \sgn(l)[ 1^{|l|}; 0^{d-|l|} ] + a\frac{\bf
j}{n} \; \; \vline \;\;  {l + a \equiv 1}\!\!\!\mod 2,\: 0 \le l n + ad < d,\: 
-\frac{d}{2} \le l < \frac{d}{2} \right\rbrace \cup \lbrace
\frac{\bf j}{n} \rbrace.
\end{equation}
Thus, we have shown that
\medskip

(A) \emph{For any point of $\Lambda_1$, minimal with respect to $\phi_{\alpha,\beta}$, there is a permutation of the components that turns it into a member of $M$ given by (\ref{M_equation}).}
 
We define a mapping from ${\mathbb R}^d$ to ${\mathbb R}^2$ by
$m(\x)=(\phi_1[\x],\phi_2[\x])$. We will call the image of $M$
the $\phi$-{\em diagram}. Lines parallel to the $\phi_1$-axis in ${\mathbb R}^2$ will
be called horizontal. Next, we will show that
\medskip

(B) \emph{If $\z_1,\z_2\in M$, then 
$m(\z_1)$ and $m(\z_2)$ on the $\phi$-diagram belong to the same
horizontal line if and only if $\z_1 \in \pm\z_2$.}
%{Furthermore, if $\z_1=\z$
%and $\z_2=-\z$  both belong to $M$, then $\z\cdot {\bf \j} = 0$.}

To prove this, take $\z= \sgn(l)[ 1^{|l|}; 0^{d-|l|} ] + a\frac{\bf \j}{n} \in
M$. If $l\ne 0$, then the condition $0 \le l n + a d < d$  uniquely determines
$a$ for a given value of $l$, and if $l=0$, then $a=1$. Therefore $l$
uniquely defines $a$.

The function $l\to |l| - \frac{l^2}{d}$ is even and increasing on
$[0,\frac{d}{2}]$, so two different points of the $\phi$-diagram may belong to the
same horizontal line if and only if their preimages are $\z_1=[ 1^{|l|},
0^{d-|l|} ] + a_1\frac{\bf \j}{n}$ and $\z_2=-[ 1^{|l|}, 0^{d-|l|} ] +
a_2\frac{\bf \j}{n}$ for some $l$,$a_1$,$a_2$. If $l\ne 0$, then $0 \le lk
+ a_1d < d$, $0 \le -lk + a_2d < d$. Adding these inequalities, we get $0
\le (a_1 + a_2)d < 2d$, so $a_1 + a_2 \in\lbrace 0,1 \rbrace$. If
$a_1+a_2=0$, then $\z_1=-\z_2$ so $m(\z_1) =
m(\z_2)$. If $a_1+a_2=1$, then the numbers $l+a_1$ and $-l +
a_2$ have different parity which contradicts conditions $l+a_1 \equiv -l +
a_2 \equiv 1 \mod 2$. If $l = 0$, then $a_1 = a_2 = 1$. This proves claim (B).

\medskip

(C) The proof of Delaunay property is based on that each line $\alpha x_1 + \beta x_2 =
1$, where $\alpha,\beta > 0$, which contains an edge of the convex hull
of the $\phi$-diagram, gives rise to a quadratic form $\phi[\x]=\alpha \phi_1[\x]
+ \beta \phi_2[\x]$ such that the ellipsoid $\phi[\x]\le1$ passes through
points of the affine lattice $\Lambda_1$ and does not have any points of $\Lambda_1$ in its interior.

We are going to show that the following set is the vertex set for a centrally-symmetric Delaunay polytope in $\Lambda_1$ with center at the origin of $\Lambda$.
\begin{equation} 
  \bV_{s,2k}^d \cup \bV_{s+1,2k}^d \cup -\bV_{s,2k}^d \cup -\bV_{s+1,2k}^d
\end{equation}
A point from $\bV_{s,2k}^d$ will be denoted by $\vv_{s,2k}$.
\medskip

(D) First consider the case of $d=7$, $2k=4$, $s = 1$. The $\phi$-diagram for these
values of parameters is shown in Figure 1 (left). We see that the
line $\frac{3}{7}x_1 + \frac{2}{3}x_2 = 1$ passes through points
$m(\vv_{1,4})=m([1,0^6])$ and
$m(\vv_{2,4})=m([{-1}^2,0^5] + \frac{\mathbf  j}{3})$, and
all other points of the $\phi$-diagram are contained in the open half-plane
$\frac{3}{7}x_1 + \frac{2}{3}x_2 > 1$. This means that polytope
$P_{1,4}^7$ is a Delaunay polytope with respect to quadratic form
$\frac{3}{7}\phi_1[\x] + \frac{2}{3}\phi_2[\x] = \frac{2}{3}\x\cdot \x +
\frac{1}{3}({\bf j}\cdot \x)^2$.

\medskip

(E) Now suppose that $d\ge 8$, $2 \le k $, and $1
\le s \le \frac{d-k-1}{2k}$. Suppose $\z = \sgn(l)[ 1^{|l|}, 0^{d-|l|} ] +
a\frac{\bf j}{n}\in M$ and $0 \le \phi_2[\z]\le \frac{d}{2k} -
\frac{d}{4k^2}$ (or, equivalently, $|l| \le \frac{d}{2k}$). We want to prove the following two statements (E1) and (E2).

(E1) \emph{If $|l| = \frac{d}{2k}$ (clearly, in this case $\frac{d}{2k}$ is an integer), then}
\begin{equation}
\pm\z = \pm\left([ 1^\frac{d}{2k}, 0^{d-\frac{d}{2k}} ] -
(\frac{d}{2k}-1)\frac{\bf j}{n}\right).
\end{equation}
A direct check shows these points belong to set $M$. Then we use claim (B)
that there are at most two points that can be written as
\[[1^{|l|}, 0^{d-|l|} ] + a\frac{\bf j}{n} \in M\] for any fixed value of $|l|$,
and there are excatly two if and only if $lk + ad = 0$, in which case the points
are $\pm\z$.

(E2) \emph{If $0\le l <
\frac{d}{2k}$, then }
\begin{equation}
\z = [ 1^{|l|}, 0^{d-|l|} ] - (l-1)\frac{\bf j}{n}.
\end{equation}
Recall that $\vv_{s,2k}$ stands for $[ 1^s, 0^{d-s} ] - (s-1)\frac{\bf j}{n}$. We have
\begin{equation}
%\begin{split}
\phi_1[\vv_{s,2k}] = {\left(s + (1-s)\frac{d}{n}\right)}^2, \quad
\phi_2[\vv_{s,2k}] = s - \frac{s^2}{d}. \\
%\end{split}
\end{equation}
All points $m(\vv_{s,2k})$ therefore belong to the parabola
\begin{equation}
\label{curve-equation} t \mapsto \left( {\left(t + (1-t)\frac{d}{n}\right)}^2,
t - \frac{t^2}{d} \right)
\end{equation}
which touches the vertical axis in the point $(0,\frac{d}{2k} -
\frac{d}{4k^2})$ when $t = \frac{d}{2k}$. The parabola is shown in dash on Figure 1. The portion of the parabola for $0 \le t \le
\frac{d}{2k}$ is the graph of a convex function.

Therefore, for each $1 \le s \le \frac{d}{2k}$ we can find a line with
equation $\alpha x_1 + \beta x_2 = 1$ which passes through points
$m(\vv_{s,2k})$, $m(\vv_{s+1,2k})$, supports the convex
hull of the $\phi$-diagram and does not contain points $m(\z)$ for
$\z\in M \setminus \{\pm \vv_{s,2k} \cup \pm \vv_{s+1,2k}\}$. Quadratic
form $\phi[\x] = \alpha \phi_1[\x] + \beta \phi_2[\x]$ defines an empty
ellipsoid centered at $\mathbf{0}$, which contains the vertices of $P^d_{s,2k}$ on its
boundary, and does not contain any other points of $\Lambda_1$.

We have proven that $P^d_{s,2k}$ is a Delaunay polytope in affine lattice
$\Lambda_1$ with respect to quadratic form $\phi=\phi^d_{s,2k}$
for $d\ge 7$. Explicit formula (\ref{formula_for_phi}) for $\phi^d_{s,2k}$
is  established by a direct calculation. 

An example of the $\phi$-diagram for $d = 19$, $2k=6$ is shown in Figure
2 (right-hand image). Note that not all points of $m(M_{rec})$ belong to
the parabola (\ref{curve-equation}).

\begin{figure}
\begin{center}
\resizebox{!}{290pt}{\includegraphics[clip=false,keepaspectratio=false]{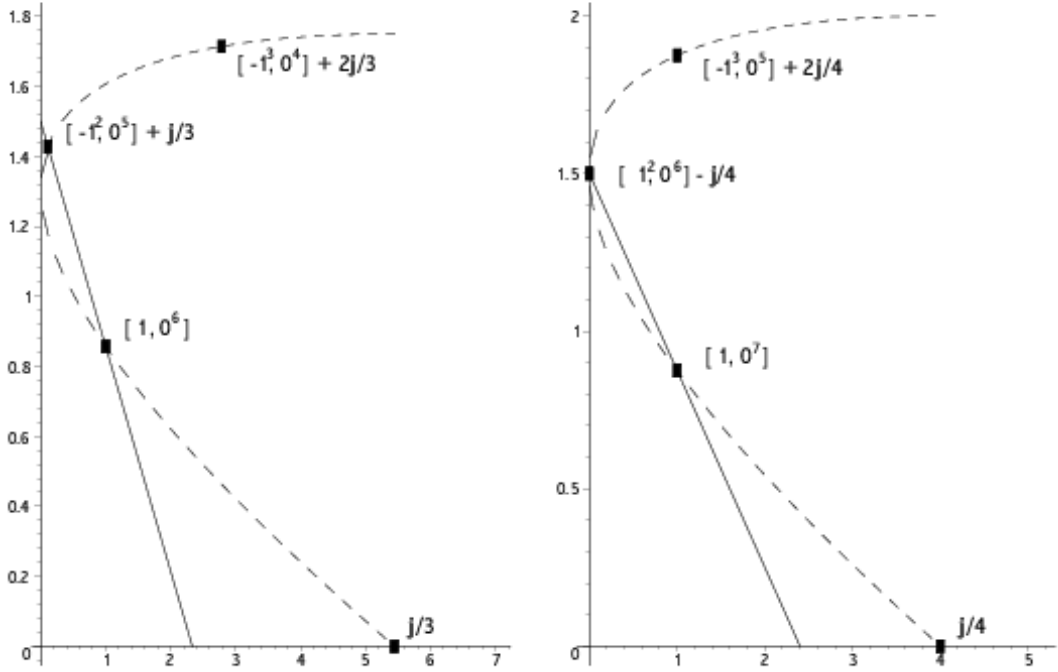}}
\caption{$\phi$-diagrams for $d=7$ (left) and $d=8$ (right); $2k=4$}
\end{center}
\end{figure}

\begin{figure}
\begin{center}
\resizebox{!}{280pt}{\includegraphics[clip=false,keepaspectratio=false]{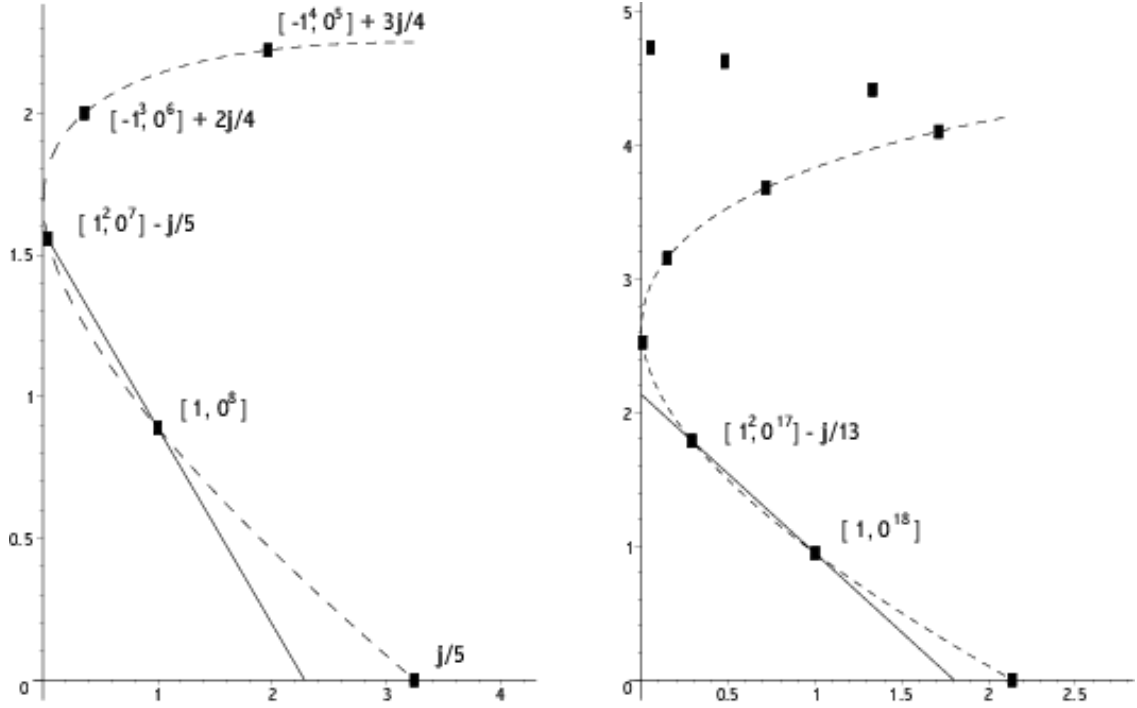}}
\caption{$\phi$-diagrams for $d=9$ (left) and $d=19$ (right); $2k=4$}
\end{center}
\end{figure}

\subsection{\large Perfection Property Part of Main Theorem}
\begin{theorem}
\label{perfectness-property-theorem} Let $d$, $k$,  $n=d-k$, $s$, $\alpha\phi_1+\beta\phi_2$ be as in Theorem
\ref{main_theorem}. Then, there is at most one 
pair consisting of a quadratic form $\phi=\alpha\phi_1+\beta\phi_2$ and $\cc \in {\mathbb R}^d$ such that the quadric $\phi[\x-\cc] = 1$ circumscribes the
polytope $P^d_{s,2k}$ .
\end{theorem}

\begin{proof}
 Since $\dim P^d_{s,2k}=d$ and  $P^d_{s,2k}$ has $\mathbf{0}$ as the center of
symmetry, $\mathbf{0}$ is the center of symmetry of the circumscribing quadric. Thus, $\cc=\mathbf{0}$. Therefore, if
$\phi$ and $\psi$ are two quadratic forms such that the corresponding quadrics 
circumscribe $P^d_{s,2k}$, then  $(\phi -
\psi)|_{\bV_{s,2k}^d \cup \bV_{s+1,2k}^d} = 0$. Let us consider an arbitrary
quadratic form $\xi$ such that $\xi|_{\bV^d_{s,2k} \cup \bV^d_{s+1,2k}} = 0$ and
prove that $\xi = 0$: this will prove the theorem.

We will use the following symmetrization technique. Let $G$ be a
subgroup of the group $S_d$ of permutations on the $d$ coordinates.  If $\tau \in G$, then we write $\tau:i \mapsto i\tau $.  The symmetrization of the form $\xi$ by $G$ is
defined as
\begin{equation}
\Sym_G\xi[\x] = \sum_{\tau\in G} \xi(x_{1\tau },\dots,x_{d\tau})=\sum_{\tau\in G} \xi[\x \tau]
\end{equation}
Since the polytope $P^d_{s,2k} \subset \R^d$ is invariant under permutations of the coordinates  of ${\mathbb
R}^d$, all components in the above
sum are $0$ on $\vertex P^d_{s,2k}$. Therefore
$\xi|_{\bV^d_{s,2k} \cup \bV^d_{s+1,2k}} = 0$.

First we prove that if a form $\xi$ satisfies the condition $\xi|_{\bV^d_{s,2k}
\cup \bV^d_{s+1,2k}} = 0$, then the sum of the diagonal elements of $f$, and
the sum of off-diagonal elements are 0. We symmetrize $f$ by the group
$G=S_d$ and get a form $\Sym_G\xi$ with diagonal coefficients proportional to the
sum of diagonal coefficients of $\xi$, and non-diagonal coefficients
proportional to the sum of non-diagonal coefficients of $f$. Therefore $\Sym_G\xi[\x]$
can be written as $\alpha \phi_1[\x] + \beta \phi_2[\x]$. Suppose
that $\alpha$ and $\beta$ are not both equal to $0$, i.e. $\alpha^2+\beta^2>0$. Since $\Sym_G\xi|_{\bV^d_{s,2k}
\cup \bV^d_{s+1,2k}} = 0$, the line $\alpha x_1 + \beta x_2 = 0$ passes
 through the points $m(\bV^d_{s,2k})$ and $m(\bV^d_{s+1,2k})$,
where $m(\x)=(\phi_1[\x],\phi_2[\x])$. This is impossible because
the line that passes through these points is uniquely defined and does not
contain $0$. This contradiction proves our claim.

Next we prove that all diagonal elements of the form $\xi$ are equal to $0$.
Suppose that one of them is nonzero. Without a limitation of generality we
may assume that $t=f_{11} \ne 0$. Let 
$St(1)$ be the subgroup of permutations which leave $1$ fixed. The matrix of
$\Sym_{St(1)}\xi$ is
\begin{equation}
\left[
\begin{matrix}
t     & \beta  & \beta  & .       & .      & \beta \cr

\beta & \alpha & \delta & .       & .      & \delta \cr

\beta & \delta  & \alpha & \delta & .      & \delta \cr

\dots \cr

\beta & \delta  & .      & \delta & \alpha & \delta \cr

\beta & \delta  & .      & .      & \delta & \alpha

\end{matrix}
\right]
\end{equation}
Consider the following vectors: \[\vv_1 = [1^s;0^{d-s}] - (s-1)\frac{\bf
j}{n} \in \bV^d_{s,2k}\quad \vv_2 = [1^{s+1};0^{d-s-1}] - s\frac{\bf j}{n} \in
\bV^d_{s+1,2k}.\] The conditions $\Sym_{St(1)}\xi(\vv_1)=0$, $\Sym_{St(1)}\xi(\vv_2)=0$ yield

\begin{multline}
t + 2(s - 1)\beta  + (s - 1)\alpha  + (s - 1)(s -2)\delta - \cr  \phantom{.} \\
\frac {2(s - 1)(t + (s + d - 2)\beta  + (s - 1)\alpha  +(s - 1)(d - 2)\delta )}{n}   \phantom{.} + \cr \\
  \frac{(s - 1)^{2}(t + 2(d- 1)\beta  + (d - 1)\alpha  + (d - 1)(d - 2)\delta )}{n^{2}} = 0,
\end{multline}

\bigskip
\bigskip

\begin{multline}
t + 2s\beta  + s\alpha  + s(s - 1)\delta -    \phantom{.}   \frac {2s(t+ (s + d - 1)\beta  + s\alpha + s(d-2)\delta )}{n} + \cr   \phantom{.} \\
\frac {s^{2} (t + 2(d - 1)\beta  + (d - 1)\alpha  + (d - 1)(d - 2)\delta
)}{n^{2}} = 0.
\end{multline}

We also have the conditions that the sums of the diagonal and the
off-diagonal elements are equal to $0$:
\begin{equation}
%\begin{split}
t + (d-1)\alpha = 0,   (d-1)(d-2)\delta + 2(d-1)\beta=0.
%\end{split}
\end{equation}
The determinant of the system of these $4$ equations in variables
$\alpha$, $\beta$, $\delta$ and $t$ is equal to
\begin{equation}
\frac{2}{n}(d - 1)(s - d + 1)(s - d)(d - 2 - n)
\end{equation}
and it is not equal to $0$ for the specified values of $d$, $n$ and $s$.
Hence $t=0$.

We have shown that the diagonal coefficients of $\xi$ are all zero. Next we prove that all off-diagonal elements of $\xi$ are also equal to $0$.
Without loss of generality we can assume that $\xi_{12}\ne 0$. Let
$\Sym_{St(1)}\xi$ be the symmetrization of $\xi$ by the group of permutations which map
the set $\lbrace 1,2 \rbrace$ onto itself. The matrix of $\Sym_{St(1)}\xi$ is then
\begin{equation}
\left[
\begin{matrix}
0     & \alpha  & \beta  & \dots       & \dots      & \beta \cr

\alpha & 0 & \beta & \dots       & \dots      & \beta \cr

\beta & \beta  & 0 & \delta & \dots      & \delta \cr

\dots & \dots  & \delta & \dots & \dots      & \dots \cr

\dots  & \dots  & \dots & \dots & \dots      & \dots \cr

\beta & \beta  & \dots      & \dots &        0            &        \delta \cr

\beta & \beta  & \delta      & \dots      & \delta & 0

\end{matrix}
\right].
\end{equation}
where $\alpha\ne 0$. We consider the following vectors: 
\[
\uu_1 =
[1^{s+1};0^{d-s-1}] - s\frac{\bf j}{n} \in \bV^d_{s+1,2k}\quad \uu_2 =
[0^{d-s-1};1^{s+1}] - s\frac{\bf j}{n} \in \bV^d_{s+1,2k}
 \]
The equations
$\Sym_{St(1)}\xi[\uu_1]= 0$, $\Sym_{St(1)}\xi[\uu_2] = 0$ yield

\begin{multline}
4(s-1)\beta + 2\alpha + (s-1)(s-2)\delta - 2\frac{s}{n}(2\alpha + 2(s-1 +
d-2)\beta + (s-1)(d-3)\delta)   + \cr \frac {(s - 1)^{2}(2\alpha  + 4(d -
2)\beta + (d - 2)(d - 3)\delta )}{n^2} = 0,
\end{multline}

\begin{equation}
%\begin{split}
(s + 1)s\delta  - \frac{2(s+1)s(2\beta  + (d - 3)\delta )}{n}   + \frac
{s^{2}(2\alpha  + 4(d - 2)\beta  + (d - 2)(d - 3)\delta )}{n^2} = 0.
%\end{split}
\end{equation}
We also know that the sum of off-diagonal elements of the matrix of $\Sym_{St(1)}\xi$ is
equal to $0$:
\begin{equation}
\label{off_diagonal_1_equation} 2\alpha  + 4(d - 2)\beta  + (d - 2)(d -
3)\delta = 0
\end{equation}
so, with some simplifications, the previous two equations can be rewritten
as

\begin{multline}\label{off_diagonal_2_equation}
4(s-1)\beta + 2\alpha + (s-1)(s-2)\delta   - 2\frac{s}{n}(2\alpha + 2(s +
d - 3)\beta + (s-1)(d-3)\delta) = 0,\cr \delta  - \frac{2(2\beta  + (d -
3)\delta )}{n} = 0.
\end{multline}

The systems of equations (\ref{off_diagonal_1_equation}) and
(\ref{off_diagonal_2_equation}) in variables $\alpha$, $\beta$ and
$\delta$ has determinant $8(d-2-n)(-d+s+1)$ which is not equal to $0$ for
the specified parameters $d$, $n$ and $s$ which proves that the system has
only zero solution. In particular, $\alpha = 0$. We have proven that all
off-diagonal elements of the matrix of form $\xi$ are equal to $0$.
\end{proof}

\par \medskip

Robert M. Erdahl
\par \url{erdahl@mast.queensu.ca}
\par Mathematics \& Statistics
\par Queen's University
\par Kingston, ON   K7L 3N6 Canada
\par \medskip
Andrei Ordine
\par Scotiabank Analytic Group
\par Toronto, ON Canada
\par \medskip
Konstantin Rybnikov
\par \url{Konstantin_Rybnikov@uml.edu}
\par University of Massachusetts at Lowell
\par Mathematical Sciences
\par Lowell, MA 01854 USA

\end{document}